\documentclass[onecolumn]{article}
\usepackage{amsmath,amssymb,amsthm}
\usepackage{graphicx} % Standard graphics package
\usepackage{fullpage}

\theoremstyle{definition}
\newcommand{\set}[1]{\{#1\}}

\newcommand{\real}{\mathbb{R}}
\newcommand{\nn}{\mathbb{N}}
\newcommand{\goesto}{\rightarrow} % sequences 
\title{
  Deriving Fa\`{a} di Bruno's formula for the derivative of a composite function via %generating functions for 
compositions of integers
} 
\author{Steffen Eger}
\date{}

\begin{document}
\maketitle
\begin{abstract}
We give yet another proof for Fa\`{a} di Bruno's formula for higher derivatives of composite functions. Our proof technique relies on reinterpreting the composition of two power series as the generating function for weighted integer compositions, for which a Fa\`{a} di Bruno-like formula is quite naturally established. 
\end{abstract}

\section{Introduction}
According to Fa\`{a} di Bruno's formula, the $n$th derivative of a composite function $G\circ F$ is given by
\begin{align}\label{eq:faa}
    \frac{d^n}{dx^n}G(F(x))=\sum\frac{n!}{b_1!\cdots b_n!}G^{(r)}(F(x))\prod_{i=1}^n\Bigl(\frac{F^{(i)}(x)}{i!}\Bigr)^{b_i},
  \end{align}
  where the sum ranges over all different solutions in nonnegative
  integers $b_1,\ldots,b_n$ of $b_1+2b_2+\cdots+nb_n=n$ and where $r$ is
  defined as $r=b_1+\cdots+b_n$. For example, for $n=3$, the three solutions for $(b_1,b_2,b_3)$ are $(0,0,1)$, $(1,1,0)$ and $(3,0,0)$, which correctly yields
\begin{align*}
  G'(F(x))\cdot F'''(x) + 3G''(F(x))\cdot F'(x)F''(x)+G'''(F(x))\cdot (F'(x))^3
\end{align*}
as third derivative of $G\circ F$. 

  Many proofs of formula \eqref{eq:faa} have been given, both based on combinatorial arguments --- such as via Bell polynomials \cite{comtet} or set partitions --- as well as on analytical; the latter, for example, using Taylor's theorem \cite{johnson}. %Actually, the formula is not very difficult to prove inductively or by application of the multinomial theorem to the composition of generating functions %\cite{flajolet}. 
%Combinatorial arguments, for instance based on the Bell polynomials or on set partitions \cite{}, are not very difficult.  
Roman \cite{roman} gives a proof based on the umbral calculus. Johnson \cite{johnson} summarizes the historical discoveries and re-discoveries of the formula as well as a variety of different proof techniques. 

%Still, 

Herein, 
we give (yet) another proof of the formula, one that is based on the combinatorics of integer compositions and a particular interpretation of the composition of power series. 
The essence of our derivation is as follows: First, we consider the number $C_{f,g}(n)$ of (doubly weighted) integer compositions of the positive integer $n$, for which we derive a closed-form formula; this requires some notation and introduction of terminology, but the derivation and combinatorial interpretation of the formula is quite intuitive. Then, 
for two arbitrary power series $G(x)=\sum_{n\ge 0} g_nx^n$  and $F(x)=\sum_{n\ge 0}f_nx^n$, we argue that $G\circ F$ has a natural interpretation of 
%we consider 
denoting 
the \emph{generating function}
\begin{align*}
   C(x)=\sum_{n\ge 0} C_{f,g}(n)x^n
\end{align*}
for $C_{f,g}(n)$. Hence, $\frac{1}{n!}\frac{d^n}{dx^n}C(0)=C_{f,g}(n)$. This yields formula \eqref{eq:faa} for $x=0$, but we argue that it is clear that the formula must indeed hold for any $x$. 

Two remarks are in order: first, as indicated, our derivation does not apply to arbitrary functions $F$ and $G$, but only to power series. 
While this may be considered a restriction, many interesting functions can indeed be represented as power series (those functions even have a name, real analytical functions). 
We also remark that, throughout, we ignore matters of convergence and treat all series as \emph{formal} and assume that functions have sufficiently many derivatives. 
Finally, while we think that many derivations of Fa\`{a} di Bruno's formula given in the literature are similar to the one we outline, we believe the particular approach that we suggest, based on integer compositions and a reinterpretation of the composition of power series, to be novel.\footnote{Technically, the approach most
similar to our own appears to be the one due to Flanders \cite{flanders}, which is, however, conceptually substantially different from our own.}
%The proof most similar to our is that given by Flanders \cite{flanders} who actually, if perceived of in the right way, performs that same derivation as we do; still, %our approach is somehwat more principled and sheds another light on Fa\`{a} di Bruno's formula in that it relates it to the integer composition problem.} 
%Finally, we remark that our proof is, if perceived of in the right way, identical to that outlined in Flanders \cite{flanders} (whose proof we found only after writing %this paper\footnote{And who makes the same power series assumptions as we do.}), who directly  Still, our approach is somewhat more principled and sheds %another light on Fa\`{a} di Bruno's formula in that it relates it to the integer composition problem. 

\section{Integer compositions and partitions}
An \emph{integer composition} of a positive integer $n$ is a tuple of positive integers $(\pi_1,\ldots,\pi_k)$, typically called \emph{parts}, whose sum is $n$. For example, the eight integer compositions of $n=4$ are 
\begin{align*}
(4),(1,3),(3,1),(2,2),(1,1,2),(1,2,1),(2,1,1),(1,1,1,1).
\end{align*}
 An \emph{integer partition} of $n$ is a tuple of positive integers $(\pi_1,\ldots,\pi_k)$ whose sum is $n$ and such that $\pi_1\ge \pi_2\ge \cdots\ge \pi_k$. For instance, there are (only) five integer partitions of $n=4$, namely 
\begin{align*}
(4),(3,1),(2,2),(2,1,1),(1,1,1,1).
\end{align*}
 Both integer compositions and partitions are well-studied objects in combinatorics \cite{andrews,heubach}. Instead of considering ordinary partitions and compositions as defined, we may consider \emph{weighted} compositions \cite{abramson,eger} and partitions, where each part value $\pi_i\in\nn=\set{1,2,3,\ldots}$ may have attributed with it a weight $f(\pi_i)\in\real$, where $\real$ denotes the set of real numbers. If weights are nonnegative integers, they may be interpreted as colors. For instance, when $f(3)=2$ and $f(1)=f(2)=f(4)=f(5)=\cdots=1$, then there are ten $f$-weighted compositions and six $f$-weighted partitions of $n=4$. These are
\begin{align*}
(4),(1,3),(1,3^*),(3,1),(3^*,1),(2,2),(1,1,2),(1,2,1),(2,1,1),(1,1,1,1)
\end{align*}
and
\begin{align*}
(4),(3,1),(3^*,1),(2,2),(2,1,1),(1,1,1,1),
\end{align*}
respectively, where we use a star ($*$) to differentiate between the two colors of part value $3$. When weights are nonintegral real numbers, they may simply be interpreted as ordinary `weights' --- possibly as probabilities if the range of $f$ is the unit interval $[0,1]$. 

Let us note that integer partitions of an integer $n$ admit an alternative, equivalent representation. Instead of writing a partition of $n$ as a tuple $(\pi_1,\ldots,\pi_k)$ with $\pi_1\ge\cdots\ge \pi_k$, we may represent it as a tuple $(k_1,\ldots,k_n)$, with $0\le k_i\le n$, for all $i=1,\ldots,n$, whereby $k_i$ denotes the \emph{multiplicity} of (type) $i\in\set{1,2,\ldots,n}$ in the composition of $n$. For instance, the above five integer partitions of $n=4$ may be represented as
\begin{align*}
(0,0,0,1),(1,0,1,0),(0,2,0,0),(2,1,0,0),(4,0,0,0).
\end{align*}
Obviously, each such tuple $(k_1,\ldots,k_n)$ must satisfy $1\cdot k_1+2\cdot k_2+\cdots+n\cdot k_n=n$. Assuming that the weighting function $f$ takes on only integral values, for the moment, how many $f$-weighted integer partitions of $n$ are there? Apparently, this number is given by
\begin{align}\label{eq:partition1}
  \sum_{k_1+2k_2+\cdots+nk_n=n}f(1)^{k_1}\dots f(n)^{k_n},
\end{align}
since the solutions, in positive numbers, of $k_1+2k_2+\cdots+nk_n=n$ are precisely the integer partitions of $n$ and the product $f(1)^{k_1}\cdots f(n)^{k_n}$ assigns the different colors to a given partition $(k_1,\ldots,k_n)$. How many $f$-weighted integer compositions of $n$ are there? Note that, in the representation $(k_1,\ldots,k_n)$ of a partition, $k_1$ denotes the number of `type' $1$, $k_2$ denotes the number of `type' $2$, ..., and $k_n$ denotes the number of `type' $n$ used in the partition of $n$. Since compositions are ordered partitions, for compositions, we need to distribute the $k_1$ types $1$, ..., $k_n$ types $n$ in a sequence of length $(k_1+\cdots+k_n)$. Therefore, the number of $f$-weighted integer compositions is simply:
\begin{align}\label{eq:composition1}
   \sum_{k_1+2k_2+\cdots+nk_n=n}\binom{k_1+\cdots+k_n}{k_1,\ldots,k_n}f(1)^{k_1}\dots f(n)^{k_n},
\end{align}
where $\binom{r}{k_1,\ldots,k_n}=\frac{r!}{k_1!\cdots k_n!}$ (for $r=k_1+\cdots+k_n$) denote the \emph{multinomial coefficients}. 

Finally, let us assume that integer partitions/compositions with a \emph{given, fixed number $k$ of parts} are (additionally) weighted (e.g., colored) by $g(k)$, for a   weighting function $g:\nn\goesto\real$. For instance, we might double count the $f$-weighted partitions/compositions with exactly $k_1+\ldots+k_n=4$ parts (or assign them higher/lower probability). Then, the number of $f$-weighted integer compositions of $n$ where parts are $g$-weighted is simply given by
\begin{align}\label{eq:combinatorial}
C_{f,g}(n)=\sum_{k_1+2k_2+\cdots+nk_n=n}\binom{k_1+\cdots+k_n}{k_1,\ldots,k_n}g(k_1+\cdots+ k_n)\prod_{i=1}^n f(i)^{k_i}.
\end{align}
 If $f$ and/or $g$ take on nonintegral values, \eqref{eq:partition1} and \eqref{eq:composition1} denote the total weight of all $f$-weighted integer partitions/compositions, and \eqref{eq:combinatorial} denotes the total weight of all $f$-weighted integer compositions of $n$ where parts are $g$-weighted. Henceforth, for brevity, we also call such compositions simply $(f,g)$-weighted.

\section{Derivation of Fa\`{a} di Bruno's formula}
We assume that $F(x)$ and $G(x)$ are the power series
\begin{align*}
  F(x) &= f_0+f_1x^1+f_2x^2+\ldots = \sum_{n\ge 0} f_nx^n,\\
  G(x) &= g_0+g_1x^1+g_2x^2+\ldots = \sum_{n\ge 0}g_nx^n,
\end{align*}
for some real coefficients $f_0,f_1,f_2,\ldots$ and $g_0,g_1,g_2,\ldots$. 
In the remainder, for ease of interpretation, we speak of the $f_n$ and $g_n$ values as if they were nonnegative and  integral, but keep in mind that they may be arbitrary real numbers. 

We interpret $F$ and $G$ as follows. The function $F$ is the generating function for the number of $f$-weighted integer compositions of $n$ with exactly one part, whereby $f(n)=f_n$. In fact, the coefficient $f_n$ of $x^n$ of $F(x)$ gives the number of $f$-weighted integer compositions of $n$ with exactly one part. We assume that $f_0=0$  (that is, integer compositions admit only positive integers).

In the context $G\circ F$, we interpret the function $G$ as follows: $G\circ F$ represents, for $G(x)=x^k$, the generating function for the number of $f$-weighted integer compositions with exactly $k$ parts; 
%in fact, this interpretation is natural since $F(x)^k$ is the 
for $G(x)=a_kx^k$, it represents the generating function for the  number of $f$-weighted integer compositions with exactly $k$ parts, where $f$-weighted compositions with $k$ parts are weighted by the factor $a_k$; and, finally, for $G(x)=x^j+x^k$, it represents the generating function for the number of $f$-weighted integer compositions with either $j$ or $k$ parts (union). This interpretation of $G$, in the context $G\circ F$, is a natural interpretation, since, for example, the coefficients of $x^n$ of $(F(x))^2$ have the form $\sum_{i=0}^nf_{n-i}f_i$, and all combinations of the number of $f$-weighted compositions of $n-i$ with one part and the number of $f$-weighted compositions of $i$ with one part yield the number of $f$-weighted compositions of $(n-i)+i=n$ with two parts.\footnote{This is in fact a critical point of our proof; if we interpreted $F(x)$ as the generating function for other combinatorial objects, such as integer partitions, then $G(x)=x^k$, in the context $G\circ F$, could not have the same interpretation as the one we have outlined.} Then, the interpretation of $(F(x))^k$ follows inductively. Similarly, if $(F(x))^k$ denotes the generating function for the number of $f$-weighted compositions of $n$ with exactly $k$ parts and $(F(x))^j$ denotes the analogous generating function for $j$ parts, then their sum obviously denotes the generating function for $k$ or $j$ parts. 

Hence, to summarize, we interpret $G\circ F$ as the generating function for the number of $f$-weighted integer compositions with arbitrary number of parts (recall that the `$+$' mean union over number of parts) and where compositions with $k$ parts are weighted by $g(k)=g_k$. Then, by virtue of the definition of generating functions, we know that $\frac{1}{n!}\frac{d^n}{dx^n}(G\circ F)(0)$ gives the number of $(f,g)$-weighted integer compositions of $n$.
%, where part number sizes $k$ are weighted by $g$. 
This is the number $C_{f,g}(n)$, whence by formula \eqref{eq:combinatorial}, we know that
\begin{align*}
\frac{1}{n!}\frac{d^n}{dx^n}(G\circ F)(0) = C_{f,g}(n)= \sum_{k_1+2k_2+\cdots+nk_n=n}\binom{k_1+\cdots+k_n}{k_1,\ldots,k_n}g(k_1+\cdots+ k_n)\prod_{i=1}^n f(i)^{k_i},
\end{align*}
or, equivalently,
\begin{align}\label{eq:interm}
  \frac{d^n}{dx^n}(G\circ F)(0) =  \sum_{k_1+2k_2+\cdots+nk_n=n}\frac{n!}{k_1!\cdots k_n!}r!g(r)\prod_{i=1}^n f(i)^{k_i},
\end{align}
where we write $r=k_1+\cdots+k_n$. 
We note that
\begin{align*}
   f(i) &= \frac{1}{i!}F^{(i)}(0), \quad\forall\:i=1,\ldots,n,\\
  r!g(r) &= G^{(r)}(0) = G^{(r)}(F(0)),
\end{align*}
whence we can rewrite \eqref{eq:interm} as
\begin{align}\label{eq:final}
  \frac{d^n}{dx^n}(G\circ F)(0) =  \sum_{k_1+2k_2+\cdots+nk_n=n}\frac{n!}{k_1!\cdots k_n!}G^{(r)}(F(0))\prod_{i=1}^n \bigl(\frac{F^{(i)}(0)}{i!}\bigr)^{k_i},
\end{align}
which is Fa\`{a} di Bruno's formula \eqref{eq:faa} evaluated at $x=0$. 
%Now, given our derivation of $\frac{d^n}{dx^n}(G\circ F)(0)$, we could either `guess' that $\frac{d^n}{dx^n}(G\circ F)(x)$ has the analogous form as %$\frac{d^n}{dx^n}(G\circ F)(0)$ and then verify, e.g., by induction. Alternatively, 
Now, 
from $(G\circ F)'(x)=G'(F(x))\cdot F'(x)$, and then $(G\circ F)''(x)=G''(F(x))F'(x)+G'(F(x))F''(x)$, etc., it %becomes apparent 
is clear that %formula \eqref{eq:interm} must in fact hold for any $x$, not only $x=0$. 
$\frac{d^n}{dx^n}(G\circ F)(x)$ is a sum of products of factors $G^{(j)}(F(x))$ and $F^{(m)}(x)$. It is also clear that, whatever the precise form of $\frac{d^n}{dx^n}(G\circ F)(x)$, evaluating it at $x=0$ will simply yield the same sum of products of factors $G^{(j)}(F(x))$ and $F^{(m)}(x)$, evaluated at $x=0$. Hence, \eqref{eq:final} must in fact hold for all $x$, not only for $x=0$. 

\section{Discussion}
We argued that $G\circ F$ has, for arbitrary power series $G$ and $F$ with coefficients $g_n$ and $f_n$, respectively, a natural interpretation as denoting the generating function for  $(f,g)$-weighted integer compositions, whereby $f(n)=f_n$ and $g(n)=g_n$, for whose coefficients Fa\`{a} di Bruno-like formulas quite effortlessly arise. 

%\section{Conclusion}

%denote a random
%draw %n 
% from the discrete
%distribution with probability mass function
\end{document}